\documentclass[11pt]{article}
\usepackage{fullpage}
\usepackage{epsfig}
\usepackage{amssymb}
\usepackage{latexsym}
\usepackage{url}
\begin{document}
\newenvironment{proof}[1][Proof]
               {\par \normalfont
                \trivlist
                \item[
                \itshape #1{.}]\ignorespaces
               }{\endtrivlist}
\newtheorem{theorem}{Theorem}[section]
\newtheorem{corollary}[theorem]{Corollary}
\newtheorem{lemma}[theorem]{Lemma}
\newtheorem{proposition}[theorem]{Proposition}
\newtheorem{remark}[theorem]{Remark}
\newtheorem{conjecture}[theorem]{Conjecture}
\newcommand{\ind}{1\hspace{-2.5mm}{1}}
\title{A useful relationship between epidemiology and queueing theory}
\author{Pieter Trapman$^{1,2}$ and Martin Bootsma$^{1,3}$}

\begin{center}{\Large {\bf A useful relationship between epidemiology and queueing theory}}\end{center}
\begin{center}{\Large{\bf Pieter Trapman$^{1,2}$ and Martinus Christoffel Jozef Bootsma$^{1,3,\dagger}$} }\end{center}

\noindent$^{1}$\emph{Julius Center for Health Research \& Primary Care; University Medical Center Utrecht, Heidelberglaan 100,  P.O.Box 85500, 3508 GA Utrecht, The Netherlands}\newline
\noindent$^{2}$\emph{Faculty of Sciences, Department of Mathematics,  Vrije Universiteit. De Boelelaan 1081a,
1081 HV Amsterdam, The Netherlands}\newline
\noindent$^{3}$\emph{Faculty of Science, Department of Mathematics,  Utrecht University. Budapestlaan 6, P.O.box 80010, 3584 CD Utrecht, The Netherlands}\newline
\quad\newline
\noindent$^{\dagger}$ Corresponding author Martin Bootsma\\
Faculty of Science, Department of Mathematics,  Utrecht University. Budapestlaan 6, P.O.box 80010, 3584 CD Utrecht, The Netherlands\\
e-mail: M.C.J.Bootsma@uu.nl\\
phone: +31-30-2531530\\
fax: +31-30-2518394\\

\newpage
\begin{abstract}
In this paper we establish a relation between the spread of infectious diseases and the dynamics of so called $M/G/1$ queues with processor sharing. The in epidemiology well known relation between the spread of epidemics and branching processes and the in queueing theory well known relation between $M/G/1$ queues and birth death processes will be combined to provide a framework in which results from queueing theory can be used in epidemiology and vice versa.

In particular, we consider the number of infectious individuals in a standard $SIR$ epidemic model at the moment of the first detection of the epidemic, where infectious individuals are detected at a constant per capita rate. We use a result from the literature on queueing processes to show that this number of infectious individuals is geometrically distributed.  
\end{abstract}
\section{Introduction}\label{intro}

\subsection{Epidemiological motivation for the paper}
In real-life, the knowledge that an infectious disease is spreading may change the further spread of the disease. Both deliberate intervention measures and non-mandatory behavioural changes can contribute to this phenomenon \cite{Boot07}. Examples are  the 2002-2003 SARS epidemic, during which people avoided crowded places; the nosocomial pathogen methicillin-resistant \emph{Staphylococcus aureus} (MRSA), where in the Nordic countries and The Netherlands hospitalised patients known to be colonised with MRSA are treated in single-bed isolation rooms and all patients who might have had contact with the index case, i.e., at least all patients in the same hospital ward as the detected patient, are screened for colonization \cite{Bootip}; or contagious animal diseases like Foot and Mouth Disease or Classical Swine Fever, where farms are depopulated at the moment the disease is detected \cite{Trap04}.  The distribution of the number of infectious individuals at the moment of the first detection of the epidemic is important. Not only because the disease dynamics before the first detection is not influenced by control measures and, therefore, can be a measure for the true spreading capacity of the disease but also because this distribution informs us on the amount of control required to curtail the outbreak.

 In this paper we deal with the change of the epidemic process at detection, by exploring an obvious but not so well-known relationship between epidemics and queueing theory. By using this relationship and results from queueing theory we are able to derive the distribution of the number of infectious individuals at the moment of the first detection in a broad class of epidemics in large populations. Furthermore, the approach of tackling problems from epidemiology by using queueing theory is promising itself and might lead to results on epidemics that are beyond the scope of this paper.

\subsection{The stochastic processes}
We first consider an $SIR$ (Susceptible $\to$ Infectious $\to$ Removed/Recovered) epidemic \cite{Ande00,Diek00} with detections in a fixed population of size $n$. In this model each pair of individuals contacts each other at a constant rate of $\lambda/(n-1)$. This implies that every individual makes contacts at a total rate $\lambda$, i.e., the number of pairs it belongs to  times the contact rate per pair. If a contact is made between an infectious and a susceptible individual, then the susceptible individual becomes infectious. Note that we define contacts as events at which infection will take place if one of the individuals is infectious and the other is susceptible. These contacts need not be the same as physical contacts. If the physical contact rate is $\lambda'/n$ and the probability of transmission during a physical contact is $c$, then $\lambda = c \lambda'$.

After an infectious period which is distributed as the random variable $L$, and independent of other infectious periods, an infectious individual becomes removed, which means that the individual  becomes immune  and stays so forever. During its infectious period an individual might be detected, which happens at a per capita rate of $\delta$. At the moment of a detection, the process might change. However, because we are interested in the number of infectious individuals at the moment of a detection, it is in the context of this paper not important what happens after the detection.

This process can be related to a branching (birth-death) process \cite{Jage75}, by interpreting an individual (say $j$) that is infected by another individual (say $i$) as a child of $i$. This interpretation relates the number of infectious individuals at time $t$ to the number of individuals in the branching process at time $t$. For finite population size $n$, this is not a proper branching process, because the distribution of children per individual changes over time because of depletion of susceptibles. However, for large $n$, the probability that during the initial phase of the epidemic contacts are made between individuals that are both non-susceptible is small and, therefore, the start of an epidemic can be analysed using the corresponding branching process \cite{Ande00,Ball95}. In this branching process individuals have independent random life lengths, all distributed as $L$. During their lifetime they get children at a constant rate $\lambda$. As in the epidemic process, detections take place at a per capita rate of $\delta$.

 \textbf{Remark:}  For definition purposes we only consider the $SIR$ epidemic in this paper. However, all results of this paper hold for all types of epidemics in which the number of infectious individuals, can be approximated by the number of living individuals in a branching process in which both the per capita birth rate and the per capita detection rate are constant over the lifetime. So, also the  $SIS$ (Susceptible $\to$ Infectious $\to$ Susceptible) epidemics, in which recovered individuals are immediately susceptible again, $SIRS$  (Susceptible $\to$ Infectious $\to$ Removed/Recovered $\to$ Susceptible) epidemics, in which recovered individuals become susceptible again after some (possibly random) time period or epidemics in homogeneous populations in which deaths and births are taken into account may be considered.

From queueing theory, we know that we can relate the so-called M(emoryless)/G(eneral)/1-queue with a Processor Sharing service discipline (or $M/G/1$-PS) queue \cite{Kell05,Kend51,Kend53,Kita93} to a branching process. 
The $M/G/1$ queue is defined as follows. Customers enter a queue at a constant rate $\lambda$ and they require a random amount of serving time (their workload) from a server. The workloads of the customers are independent and distributed as $L$. A single server is serving the queue.

The processor sharing service discipline may be interpreted as the limit of a round robin service discipline. In that discipline, the server serves a customer either for a time length $h$, or if the remaining amount of serving time is less than $h$, until the customer is completely served. After this time length the server moves on to the next customer. Customers that are not completely served yet stay in the queue, while completely served customers leave the queue. After serving the last customer in the queue, the server returns to the first customer that is still in the queue. The processor sharing discipline is obtained by taking the limit $h \searrow 0$.

For the $M/G/1$ queue with round robin service discipline, we may consider customers that arrive during the time that the server was serving customer $i$ as children of $i$. In this way we obtain a branching process, in which individuals can give birth to other individuals during a random period distributed as $L$ and during which it gets children at a constant rate $\lambda$. Note that time in the queueing process is not the same as time in the branching process, because in the branching process the individuals all grow older at a constant speed, while in the queueing process, customers are served one by one, where the server has a constant speed. So, if the number of customers in the row is increasing, the service per time unit for a given individual in the queue is decreasing. However, if $h \searrow 0$, the order of events (arrivals/births and departures/deaths) in the processes is the same. Another difference between the queueing process and the corresponding branching process is that for the branching process the state without individuals is absorbing, while in the queueing process customers may arrive in an empty queue, but up to the first time the queue becomes empty, the two processes are in direct correspondence.

We add a ``catastrophe'' process \cite{Stir06} to the queueing process in which catastrophes occur at a constant rate $\delta$. At the time of a catastrophe, a sudden change might happen to the queue, e.g., the queue might be emptied or a fraction of the queue might leave the queue. However, in this paper we are only interested in the number of customers in the queue at the moment of a catastrophe. Therefore, we do not need to specify what happens at the time of a catastrophe. Yet in the epidemic context we shall assume that detection has no impact at all, as this helps to keep the notation simple.

The catastrophes are incorporated in the branching process corresponding to the $M/G/1$-PS queue as follows: Consider the $M/G/1$ queue with round robin service discipline. If a catastrophe occurs at a moment customer $i$ is served, then it can be seen as a detection of individual $i$ in the branching process, which occur at a per-capita rate $\delta$.  The $M/G/1$-PS queue with catastrophes is obtained by $h \searrow 0$, where $h$ is the time a server stays with the same customer. The order of events (arrivals/births, departures/deaths and catastrophes/detections) in the branching process with detections and the $M/G/1$-PS queue is still maintained. So, the number of individuals alive at the moment of the first detection in a branching process has the same distribution as the number of customers in the $M/G/1$-PS queue at the first catastrophe. 

We use a result from Kitaev \cite{Kita93}, to show that the number of customers in a queue after an exponentially distributed time, conditioned on the queue being in its first busy period, i.e., the period in which the server was non-stop working, starting at the arrival of the first customer,  is geometrically distributed. This implies that the number of infectious individuals at the moment of the first detection for an $SIR$ epidemic in a large randomly mixing population is also geometrically distributed. So, the distribution of this number can be described by one parameter, while the process itself is described by 2 parameters and an unspecified distribution. This implies that if the only observations available, are the number of infectious individuals at the moment of the first detection in different outbreaks of the same disease, then we can only estimate 1 parameter and we cannot provide estimates for $\lambda$, $\delta$ and $L$ separately. 

\subsection{Earlier work}

The relation between $M/G/1$ queues and birth and death processes has already been discussed by Kendall in \cite{Kend51,Kend53}. In \cite{Kita93} a discussion on the relation between birth and death processes and the $M/G/1$ queues with processor sharing can be found. References to earlier work on this subject are also given in that paper. 

We are aware of only a few references to queueing theory in the epidemiological literature. In \cite{Ball95} results for the $M/M/1$ queue (the queue with exponentially distributed workloads) are used to compute the total cost of an epidemic, which is interpreted as a constant times the total number of infection days (i.e., it is interpreted as a constant times the integral of the number of infecteds over time). In \cite{Sell83} (see also \cite[p.12]{Ande00}) a construction for an $SIR$ epidemic is given, which is very close to interpreting the epidemic as an $M/G/1$ queue with processor sharing, however the relation is not made explicit there. The random time change discussed here was also applied in \cite{Wats80}, but in that paper was no explicit reference to queueing theory. 

\subsection{Outline of the paper}

In the next section, we define the $SIR$ epidemic with detections, the corresponding branching process with detections and the $M/G/1$-PS queue with catastrophes in   mathematical language, in order to show the relationship between epidemics and queues in a rigorous way. 
In Section \ref{secrela} we will provide a random time change argument to establish the relation between $SIR$ epidemics and $M/G/1$-PS queues and we show that the number of infectious individuals at the moment of first detection in a large population is geometrically distributed. We use explicit computations or refer to literature in which the results are already proven rigorously. However, if possible, we also give heuristic and intuitive arguments for the claims made, which might be helpful for further use of queueing theory in epidemiology. In particular, in Section \ref{secspec} we provide an intuitive proof for the geometric distribution of the number of infectious individuals at the moment of first detection if the infectious periods are exponentially distributed.
Some applications of the results of this paper are discussed in Section \ref{secapli}. In the final section we discuss possible extensions and some limitations of the use of queueing theory in epidemics. In particular, we discuss whether the coupling between branching processes and M/G/1 queues can be made for branching processes corresponding to epidemic models in which the infectivity/contact rate of an individual is not constant during its infectious period. 

\section{Definitions and notation}\label{secnot}

Before we formally describe the relevant processes for this paper, we give some definitions. We define $1/0 := \displaystyle\lim_{x \searrow 0} 1/x = \infty$, $1/\infty := \displaystyle\lim_{x \to \infty} 1/x=0$ and $0 \times \infty = \displaystyle\lim_{x \to \infty} 0 \times x = 0$. 
Furthermore, the maximum/supremum of an empty set of real numbers is 0, while the minimum/infimum of an empty set is $\infty$. 
The indicator function $\ind(A)$ takes value 1 if the event $A$ occurs and 0 if the event $A$ does not occur.   
For a function $f(x)$, we define $f(x-) := \displaystyle\lim_{y \nearrow x} f(y)$ and $f(x) = o(x)$ if $\displaystyle\lim_{x \searrow 0} f(x)/x = 0$. 
The natural numbers, not including $0$, are denoted by $\mathbb{N}$ and  $\mathbb{N}_0 := \mathbb{N} \cup \{0\}$.
Finally, with some abuse of notation, for all processes under consideration, we will use $\{\mathcal{F}_t;t \geq 0\}$ to denote the filtration to which the process is adapted \cite[p.~475]{Grim92}. So, all information on the process available at time $t$ is contained in $\mathcal{F}_t$.

Throughout this paper we may deviate from standard notation in literature, because we want to relate different processes, which all have different standard notations for quantities that are related in this paper. In our notation we want to make clear which quantities in the different processes are related.

\subsection{The SIR epidemic with detections}

We consider an $SIR$ epidemic in a homogeneous and randomly mixing population without demographic turnover of size $n$. Let $S^{(n)}(t)$ be the number of susceptible individuals in the population at time $t$, $I^{(n)}(t)$ the number of infectious individuals at time $t$ and $R^{(n)}(t)$ the number of removed/recovered individuals at this time. The epidemic starts with one infectious individual in a further susceptible population, i.e., $S^{(n)}(0) = n-1$, $I^{(n)}(0)=1$ and $R^{(n)}(0)=0$. Furthermore, we assume that the initial infectious individual was infected itself at time $0$.

Every pair of individuals makes contact at a constant, strictly positive, rate $\lambda/(n-1)$,  i.e., contacts between a pair of individuals are made according to a Poisson process with parameter $\lambda/(n-1)$. If a contact is made between an infectious and a susceptible individual, the susceptible individual becomes immediately infectious.

An infected individual stays infectious for a random infectious period, which is distributed as the  random variable $L$, which is almost surely (a.s.) positive, i.e., $\mathbb{P}(0<L \leq \infty)=1$. The infectious periods are independent and identically distributed (i.i.d.). After the infectious period individuals become immune and stay so forever.

We extend the standard $SIR$ epidemic model by adding a detection process. In this process, infectious individuals are detected at a constant per capita rate $\delta$. Mathematically, the detection process corresponds to a Poisson process with parameter $\delta > 0$, which is defined on all infectious periods. Let $D^{(n)}_i$ denote the random time of the $i$-th detection in the population. If the number of infectious individuals is 0 before the $i$-th detection takes place, we say $D^{(n)}_i = \infty$.  As mentioned before, in this paper we are interested in the number of infectious individuals at the time of first detection $I^{(n)}(D^{(n)}_1)$, conditioned on $D^{(n)}_1< \infty$.

\subsection{The branching process}

Let $Z(t):=Z_{\lambda,\beta(s),\delta}(t)$ be a branching process \cite{Jage75} with detections in which individuals have i.i.d.\ life lengths, distributed as $L$, with moment generating function $\beta(s) := \mathbb{E}(e^{-s L})$ and give birth at a constant per capita rate $\lambda > 0$. So new individuals are born at a total rate of  $\lambda Z(t)$. The process starts with one individual, whose age at time $0$ is $0$.  There is a detection process on top of  this process of births and deaths, which is a homogeneous Poisson process with rate $\delta > 0$, which is defined on the life times of the individuals. So, detections happen at a total non-constant rate $\delta Z(t)$. Again, detections do not influence the further course of the branching process. We only use the subscripts in $Z_{\lambda,\beta(s),\delta}(t)$ if we want to stress the dependence on these parameters.

We use $D_i$ to denote the time of the $i$-th detection in real-time in the branching process. If the branching process goes extinct before the $i$-th detection, then $D_i = \infty$. We label the individuals in the epidemic by the real-time order in which they appear in the branching process, so the ancestor in the branching process gets label 1, its first child gets label 2, then the next individual to be born (which is either a child or grand-child of the ancestor) gets label 3, and so on. Note that in general, it is not individual $i$ that is detected at time $D_i$. The random life length of individual $i$ is denoted by $L_i$ and its time of birth by $T_i$. So, $T_1=0$. We define $A_i(t)$, the age of individual $i$ at time $t$, as $A_i(t) = (t- T_i) \ind(T_i < t)$ for $t<T_i+L_i$ and $A_i(t) = \infty $ if $t \geq T_i + L_i$. So, $\frac{d}{dt}A_i(t) = \ind(0 < A_i(t) < L_i)$.

Observe that
\begin{equation}\label{branchderi1}
Z(t) = \sum_{i=1}^{\infty} (\ind(T_i \leq t) - \ind(A_i(t)=\infty))
\end{equation}
is defined in terms of events occurring no later than $t$, so $Z(t) \in \mathcal{F}_t$.

For $i \in \mathbb{N}$, the random variables $T_i$ and $D_i$ are such that both of the sequences are increasing in $i$ and
\begin{equation}\label{branchderi2}
\begin{array}{rcl}
\mathbb{P}(T_i \in (t,t+h)|\mathcal{F}_t, Z(t)=k, \displaystyle\sum_{j=1}^{\infty} \ind(T_j \leq t)= i-1) & = & \lambda k h +o(h),\\
\mathbb{P}(D_i \in (t,t+h)|\mathcal{F}_t, Z(t)=k, \displaystyle\sum_{j=1}^{\infty} \ind(D_j \leq t)= i-1) & = & \delta k h +o(h),\\
\mathbb{P}(\displaystyle\sum_{i=1}^{\infty} \ind(T_i \in (t,t+h)) + \displaystyle\sum_{i=1}^{\infty} \ind (D_i \in (t,t+h)) > 1) & = & o(h).
\end{array}
\end{equation}
If $Z(t) = 0$ and  $\displaystyle\sum_{j=1}^{\infty} \ind(T_j \leq t)= i$ (resp.\ $\displaystyle\sum_{j=1}^{\infty} \ind(D_j \leq t)= i$), then $T_k = \infty$ (resp.\ $D_k = \infty$) for $k>i$.

\subsection{The $M(\lambda)/G(\beta(s))/1$ queue with Processor Sharing and catastrophes}

We define the $M(\lambda)/G(\beta(s))/1$ queue with Processor Sharing with catastrophes or $M/G/1$-PS queue with catastrophes, $Q(t) := Q^{\mbox{PS}}(t) := Q^{\mbox{PS}}_{\lambda,\beta(s),\delta}(t)$, as follows.
Customers arrive in a queue according to a homogeneous Poisson process on the positive half line $(0,\infty)$, with rate $\lambda$. Each customer brings in a workload, which is independent of the arrival process and workloads of other customers and distributed as the random variable $L$, with moment generating function $\beta(s)$. We assume that $\mathbb{P}(0< L \leq \infty)=1$. The customers in the queue are labelled according to the order in which they arrived in the queue. The time of arrival of the $i$-th customer is denoted by $\tilde{T}_i$. Unless specified otherwise, we assume that no customers are in the queue at time $t=0$.

The total workload customer $i$ brings in is denoted by $L_i$. The amount of service time customer $i$ already received at time $t$, is denoted by $\tilde{A}_i(t)$, where we define $\tilde{A}_i(t)= 0$ for $t \leq \tilde{T}_i$. Customer $i$ leaves the queue at $\sup\{t;\tilde{A}_i(t)<L_i\}$ and from that time on (including the time of departure) $\tilde{A}_i(t)=\infty$.

One server serves the people in the queue in such a way that all customers in the queue receive the same amount of service per time unit. Thus, $\frac{d}{dt}\tilde{A}_i(t) = \ind(0 < \tilde{A}_i(t) < \infty)/Q(t)$, where $Q(t)$ is the number of customers in the queue at time $t$. 

Independently of the ordinary $M/G/1$-PS queueing process as defined above, we define a catastrophe process. Catastrophes occur according to a Poisson process on $(0,\infty)$ with rate $\delta$. The time of the $i$-th catastrophe is denoted by $\tilde{D}_i$.

Note that
\begin{equation}\label{mg1deri2}
Q(t) = \sum_{i=1}^{\infty} (\ind(\tilde{T}_i \leq t) - \ind(\tilde{A}_i(t)=\infty)).
\end{equation}
So, $Q(t)$ is (as is $Z(t)$) defined in terms of events occurring no later than $t$, i.e., $Q(t) \in \mathcal{F}_t$.
For $i \in \mathbb{N}$, the random variables $\tilde{T}_i$ and $\tilde{D}_i$ are defined such that both of the sequences are increasing in $i$ and
\begin{equation}\label{mg1deri3}
\begin{array}{rcl}
\mathbb{P}(\tilde{T}_i \in (t,t+h)|\mathcal{F}_t, \displaystyle\sum_{j=1}^{\infty} \ind(\tilde{T}_j \leq t)= i-1) & = & \lambda h +o(h),\\
\mathbb{P}(\tilde{D}_i \in (t,t+h)|\mathcal{F}_t, \displaystyle\sum_{j=1}^{\infty} \ind(\tilde{D}_j \leq t)= i-1) & = & \delta h +o(h),\\
\mathbb{P}(\displaystyle\sum_{i=1}^{\infty} \ind(\tilde{T}_i \in (t,t+h)) + \displaystyle\sum_{i=1}^{\infty} \ind(\tilde{D}_i \in (t,t+h)) > 1) & = & o(h).
\end{array}
\end{equation}
Note that, while Z(0) $\neq 0$, $Q(0)=0$.
 
\section{The relationship between the SIR-epidemic with detections and the $M/G/1$-PS queue with catastrophes}\label{secrela}

In this section we show that for the $SIR$ epidemic with infection rate $\lambda$, detection rate $\delta$ and moment generating function of the infectious period $\beta(s)=\mathbb{E}(e^{-sL})$, 
$$\lim_{n \to \infty} \mathbb{P}(I^{(n)}(D_1^{(n)})=k|D_1^{(n)}<\infty) = \mathbb{P}(Z_{\lambda,\beta(s),\delta}(D_1)=k|T_1=0,D_1 < \infty).$$
So, in the large population limit, the number of infected individuals at the time of first detection is distributed as the number of alive individuals at the time of first detection in the corresponding branching process. 

After that we prove that 
\begin{equation}\label{mainequa}
\mathbb{P}(Z_{\lambda,\beta(s),\delta}(D_1)=k|T_1=0, D_1 < \infty) = \mathbb{P}(Q^{\mbox{PS}}_{\lambda,\beta(s),\delta}(\tilde{D}_1)=k|Q^{\mbox{PS}}_{\lambda,\beta(s),\delta}(0) = 0,Q^{\mbox{PS}}_{\lambda,\beta(s),\delta}(\tilde{D}_1)>0).
\end{equation}

and finally we use this result to show that there is a $p := p(\lambda,\beta(s),\delta)$ such that for $k \in \mathbb{N}$,
\begin{equation}\label{geomresu}
\mathbb{P}(Z_{\lambda,\beta(s),\delta}(D_1)=k|T_1=0, D_1 < \infty) = p(1-p)^k.
\end{equation}

\subsection{The relationship between the epidemic process and the branching process}

Ball and Donnelly \cite{Ball95} proved that the epidemic process $I^{(n)}(t)$ for $n \in \mathbb{N}$ and the branching process $Z(t)$ can be coupled in such a way that there exists a constant $c>0$ such that 
\begin{displaymath}
\sup_{0 < t < c \log(n)} |I^{(n)}(t)-Z(t)| \to 0 \qquad \mbox{a.s.\ for $n \to \infty$.}
\end{displaymath}
Furthermore, observe that if $D^{(n)}_1< \infty$, then $I^{(n)}(t) \geq 1$ for $0<t<D^{(n)}_1$, and 
$$\mathbb{P}(D^{(n)}_1<t|D^{(n)}_1<\infty) \geq 1 - e^{-\delta t}.$$ 
This, in turn, implies that for $c>0$,
\begin{equation}\label{logbound}
\lim_{n \to \infty} \mathbb{P}(D^{(n)}_1< c \log(n) |D^{(n)}_1<\infty) =1.
\end{equation}

Combining the result by Ball and Donnelly with (\ref{logbound})  gives that for any $k \in \mathbb{N}$,
\begin{eqnarray*}
\ & \ & \lim_{n \to \infty} \mathbb{P}(I^{(n)}(D^{(n)}_1) = k|D^{(n)}_1<\infty)\\
\ & = &  \lim_{n \to \infty} \mathbb{P}(I^{(n)}(D^{(n)}_1) = k|D^{(n)}_1<c \log(n)) \mathbb{P}(D^{(n)}_1 < c \log(n)|D^{(n)}_1 < \infty)\\
\ & \ & + \lim_{n \to \infty} \mathbb{P}(I^{(n)}(D^{(n)}_1) = k|D^{(n)}_1 \geq c \log(n)) \mathbb{P}(D^{(n)}_1 \geq c \log(n)|D^{(n)}_1 < \infty)\\
\ & = &  \lim_{n \to \infty} \mathbb{P}(I^{(n)}(D^{(n)}_1) = k|D^{(n)}_1<c \log(n))\\
\ & = &  \lim_{n \to \infty} \mathbb{P}(Z(D_1) = k|D_1^{(n)}<c \log n)\\
\ & = & \mathbb{P}(Z(D_1) = k|D_1<\infty).
\end{eqnarray*}
The above results allow us to analyse the branching process $Z(t)$, instead of the $SIR$ epidemic in large populations. So, from now on we will consider $Z(t)$ instead of $I^{(n)}(t)$.

\subsection{A random time change: from branching processes to the $M/G/1$-PS queue}\label{timegeneral}

In this subsection we use a random time change to show that equation (\ref{mainequa}) holds.
Let $$\tau(t) =  \int_0^t 1/Z(t') dt'.$$ 
We see that for $i \in \mathbb{N} $ and  $\tau(t) \leq T_i$,  $A_i(\tau(t)) = 0$. Furthermore, 
$$\frac{d}{dt}A_i(\tau(t)) = \ind(0 < A_i(\tau(t)) < \infty) \frac{d\tau(t)}{dt} = \ind(0 < A_i(\tau(t)) < \infty)/Z(t),$$
where $T_1=0$ and for $i<j$, $T_i < T_j$ a.s. The random time change does not change the fact that the random variables $L_i$ are i.i.d., distributed as $L$ and independent of $T_j$ for $1 \leq j \leq i$.
Furthermore,
\begin{equation}\label{coupderi1}
Z(\tau(t)) = \sum_{i=1}^{\infty} (\ind(T_i \leq \tau(t)) - \ind(A_i(\tau(t)) = \infty)),
\end{equation}
and
\begin{equation}\label{coupderi2}
\begin{array}{rl}
\ & \mathbb{P}(T_i \in (\tau(t),\tau(t+h))|\mathcal{F}_{\tau(t)}, Z(\tau(t))=k, \sum\limits_{j=1}^{\infty} \ind(T_j \leq \tau(t))= i-1)\\
= & \lambda k (\tau(t+h)-\tau(t)) +o(h)\\
= & \lambda k (h/k +o(h)) +o(h)\\
= &\lambda h \ind(k>0) +o(h).\\
\end{array}
\end{equation}
Similarly, we deduce
\begin{equation}\label{coupderi3}
\begin{array}{rcl}
\mathbb{P}(D_i \in (\tau(t),\tau(t+h))|\mathcal{F}_{\tau(t)}, Z(t)=k, \sum\limits_{j=1}^{\infty} \ind(D_j \leq \tau(t))= i-1) & = &  \delta h \ind(k>0) +o(h),\\
\mathbb{P}(\sum\limits_{i=1}^{\infty} \ind(T_i \in (\tau(t),\tau(t+h))) + \sum\limits_{i=1}^{\infty} \ind (D_i \in (\tau(t),\tau(t+h))) > 1) & = &  o(h).
\end{array}
\end{equation}

Note that, as long as $Z(\tau(t))>0$, the description of the process $Z_{\lambda,\beta(s),\delta}(\tau(t))$ is the same as the description of $Q^{\mbox{PS}}_{\lambda,\beta(s),\delta}(t)$, with $T_i$ (resp.\ $D_i$) replaced by $\tilde{T}_i$ (resp.\ $\tilde{D}_i$). So,
\begin{equation}\label{branque1}
\mathbb{P}(Z_{\lambda,\beta(s),\delta}(D_1)=k|T_1=0,Z_{\lambda,\beta(s),\delta}(D_1)>0) = \mathbb{P}(Q^{\mbox{PS}}_{\lambda,\beta(s),\delta}(\tilde{D}_1)=k|\tilde{T}_1=0, \min_{0 \leq t \leq \tilde{D}_1} Q^{\mbox{PS}}_{\lambda,\beta(s),\delta}(t) > 0).
\end{equation}

We proceed by showing that 
$$\mathbb{P}(Q(\tilde{D}_1)=k|\tilde{T}_1=0,\min_{0 \leq t \leq \tilde{D}_1}Q(t)>0) = \mathbb{P}(Q(\tilde{D}_1)=k|Q(0)=0,Q(\tilde{D}_1)>0).$$
We use the following notation for the starting and stopping times of the busy periods of the queue: 
\begin{eqnarray*}
\sigma_1 & := &  \min\{t \geq 0;Q(t)>0\}, \\
\sigma^*_n & :=  & \min\{t>\sigma_n; Q(t)=0\}, \qquad \mbox{for $n \in \mathbb{N}$.}\\
\sigma_{n} & := & \min\{t>\sigma^*_{n-1};Q(t)>0\} \qquad \mbox{for $n \in \mathbb{N}\setminus\{1\}$.} \\
\end{eqnarray*}
Furthermore, let $N(t) := \max\{n \geq 1; \sigma^*_n \leq t \}$, be the number of times the queue becomes empty in the interval $(0,t)$.
Observe that $$\mathbb{P}(Q(t+s)=k|\sum_{i=1}^{\infty}\ind(\sigma_i=s)=1,\mathcal{F}_s) = \mathbb{P}(Q(t)=k|T_1=0),$$
and that $\tilde{D}_1$ is exponentially distributed and independent of the queue length upto time $\tilde{D}_1$. Therefore,
\begin{equation}\label{relationtwoqueueing}
\begin{array}{rcl}
\ & \ & \mathbb{P}(Q(\tilde{D}_1)=k|Q(0)=0,Q(\tilde{D}_1)>0)\\
\ & = & \sum_{n=0}^{\infty} \mathbb{P}(Q(\tilde{D}_1)=k|N(\tilde{D}_1)=n, Q(0)=0,Q(\tilde{D}_1)>0)\\
\ & \ & \mathbb{P}(N(\tilde{D}_1)=n| Q(0)=0,Q(\tilde{D}_1)>0)\\
\ & = & \sum_{n=0}^{\infty} \mathbb{P}(Q(\tilde{D}_1)=k|N(\tilde{D}_1)=n,  0 < \sigma_{n+1} < \tilde{D}_1, Q(\tilde{D}_1)>0)\\
\ & \ & \mathbb{P}(N(\tilde{D}_1)=n| Q(0)=0,Q(\tilde{D}_1)>0)\\
\ & = & \sum_{n=0}^{\infty} \mathbb{P}(Q([\tilde{D}_1-\sigma_{n+1}]+\sigma_{n+1})=k|N(\tilde{D}_1)=n,  0 =\tilde{T}_1 < \sigma_{n+1} < \tilde{D}_1, Q(\tilde{D}_1)>0)\\
\ & \ & \mathbb{P}(N(\tilde{D}_1)=n| Q(0)=0,Q(\tilde{D}_1)>0)\\
\ & = & \sum_{n=0}^{\infty} \mathbb{P}(Q(\tilde{D}_1)=k|N(\tilde{D}_1)=0,  \tilde{T}_1=0 , Q(\tilde{D}_1)>0)\\
\ & \ & \mathbb{P}(N(\tilde{D}_1)=n| Q(0)=0,Q(\tilde{D}_1)>0)\\
\ & = & \mathbb{P}(Q(\tilde{D}_1)=k| \tilde{T}_1=0, \tilde{D}_1 \leq \sigma^*_1)\\
\ & = & \mathbb{P}(Q(\tilde{D}_1)=k|\tilde{T}_1 = 0, \min_{0\leq t \leq \tilde{D}_1}Q(t)>0).
\end{array}
\end{equation}
Combining this result with equation  (\ref{branque1}) leads to equation (\ref{mainequa}).

\subsection{The distribution of $Q(\tilde{D}_1)$ in the $M/G/1$-PS queue}\label{psdist}

By \cite[eq.~(2.6)]{Kita93}, we know that for $Q(t) = Q^{\mbox{PS}}_{\lambda,\beta(s),\delta}$,
$$\int_0^{\infty} e^{-\delta t} \mathbb{E}(s^{Q(t)}|Q(0)=0) dt = \frac{1}{\delta + (1-s)\lambda(1-\pi)},$$
where $\pi$ is the smallest  root of the equation $\pi = \beta(\delta + \lambda (1-\pi)) := \mathbb{E}(e^{-(\delta +\lambda(1-\pi)){L}})$. Since $g(x) := \mathbb{E}(e^{-(\delta +\lambda(1-x)){L}})-x$, is convex and $g(0) > 0$ and $g(1) <0$, $\pi$ is the unique root of $g(x) =0$ in $[0,1]$.

We observe that 
$$\mathbb{E}(s^{Q(\tilde{D}_1)}|Q(0)=0) = \delta \int_0^{\infty} e^{-\delta t} \mathbb{E}(s^{Q(t)}|Q(0)=0) dt =  \frac{\delta}{\delta + (1-s)\lambda(1-\pi)}$$
and that if the random variable $X$ is geometrically distributed with parameter $p$, then $\mathbb{E}(s^{X-1}) = \sum_{k=0}^{\infty} p(1-p)^k s^k = p(1-(1-p)s)^{-1}$. 
By combining these observations with the fact that the probability generating function determines a distribution on the positive integers completely \cite{Grim92}, we deduce that 
\begin{equation}
\mathbb{P}(Q(\tilde{D}_1)=k|Q(0)=0,Q(\tilde{D}_1)>0) = p(1-p)^{k-1}
\end{equation}
with 
\begin{equation}\label{pandpi}
\begin{array}{rcl}
p & = & \delta/(\delta + (1-\pi) \lambda),\\
\pi & = & \{x \in [0,1]; x=\mathbb{E}(e^{-(\delta + (1-x)\lambda)L)})\},
\end{array}
\end{equation}
i.e., conditioned on $Q(0) =0$ and $Q(\tilde{D}_1)>0$, $Q(\tilde{D}_1)$ is geometrically distributed (this has already been observed in \cite{Kell05}). 

Combined with the results of the previous subsection, we come to the main result of this paper: For $k \in \mathbb{N}$
\begin{equation}
\mathbb{P}(Z_{\lambda,\beta(s),\delta}(D_1)=k|D_1<\infty) = p(1-p)^{k-1},
\end{equation}
with $p$ as above. 

\section{A special case: Markovian models}\label{secspec}

The proof of the result by Kitaev \cite{Kita93} discussed in \ref{psdist} is rigorous, but we did not succeed in finding a intuitive argument for why it should be true. However, if the queueing process has the Markov property, then we can provide an intuitive proof. This will be done in this section. 

If $L$ is exponentially distributed with parameter $\mu$, then the $SIR$ epidemic, the branching process and the queueing process under consideration are all Markovian.
The branching process then becomes a simple birth-death process \cite[p.251]{Grim92} with detections, described by the following equations: 
\begin{equation}\label{birthdeatheq}
\begin{array}{rcl}
\mathbb{P}(Z(t+h)= k+1|Z(t)=k) & = & \lambda k h +o(h),\\ 
\mathbb{P}(Z(t+h)= k-1|Z(t)=k) & = & \mu k h +o(h),\\
\mathbb{P}(D_i \in (t,t+h)|D_{i-1} \leq t < D_i, Z(t)=k) & = & \delta k h + o(h),\\
\mathbb{P}(\mbox{more than 1 event in $(t,t+h)$}|Z(t) = k) & = & o(h).
\end{array} 
\end{equation}

The $M(\lambda)/G(\frac{\mu}{\mu+s})/1$-queue with exponentially distributed workloads with parameter $\mu$, is usually referred to as the $M(\lambda)/M(\mu)/1$ queue \cite[p.420]{Grim92},\cite{Kend51}. A property of this model is that many service disciplines lead to the same process $Q(t)$, because no matter which customer is served, the served customer will leave the queue with probability $\mu h + o(h)$ during an interval of length $h$. In particular, the First-In-First-Out (FIFO) discipline (in which customers are served 1 by 1 until they are fully served in the order of arrival), the Last-In-First-Out (LIFO) discipline (in which the server always serves the customer that arrived last out of the customers still in the queue) and the Processor Sharing discipline for $M/M/1$ queues with catastrophes all lead to the equations   
\begin{equation}\label{mm1eq}
\begin{array}{rcl}
\mathbb{P}(Q(t+h)= k+1|Q(t)=k) & = & \lambda h +o(h),\\
\mathbb{P}(Q(t+h)= k-1|Q(t)=k) & = & \mu h \ind(k>0) +o(h),\\
\mathbb{P}(\tilde{D}_i \in (t,t+h)|\tilde{D}_{i-1} \leq t < \tilde{D}_i) & = & \delta h + o(h),\\
\mathbb{P}(\mbox{more than 1 event in $(t,t+h)$}) & = & o(h).
\end{array} 
\end{equation}

We continue by considering an $M/G/1$ queue with catastrophes if the server uses a LIFO discipline, $Q^{\mbox{LIFO}}(t)=Q_{\lambda,\beta(s),\delta}^{\mbox{LIFO}}(t)$. 
We deduce  that for $k \in \mathbb{N}$ and some $0 < p' <1$, 
$$\mathbb{P}(Q^{\mbox{LIFO}}(\tilde{D}_1)=k|\tilde{D}_1 < \infty) = p'(1-p')^k.$$ 
It will turn out that $p'=p$, \cite{Kell05}.
If $L$ is exponentially distributed, then this result applies also to the processor sharing discipline, because as stated above, $Q^{\mbox{LIFO}}(t)$ and $Q^{\mbox{PS}}(t)$ have the same law for $M/M/1$ queues. 

\subsection{The distribution of $Q(\tilde{D}_1)$ in the $M/G/1$-LIFO queue with catastrophes}\label{lifodist}

For the $M/G/1$-LIFO queue with catastrophes, the event $\{\tilde{T}_1=0,\displaystyle\min_{0 \leq t \leq \tilde{D}_1}Q(t)>0\}$ is the same as the event $\{\tilde{T}_1=0,\mbox{customer $1$ is still in the queue at time $\tilde{D}_1$}\}$. Furthermore, note that if a LIFO service discipline is applied and customer $i$ is not in the queue at time $\tilde{D}_1$, then none of the customers that arrived between the arrival and departure of customer $i$, will be in the queue at time $\tilde{D}_1$. So, we may ignore the arrival of such customers.
Therefore, we say that the only two ``good'' events are arrivals of new customers, that will be in the queue at time $\tilde{D}_1$ and catastrophes. The first ``good'' event after $\tilde{T}_1$ occurs during the period that customer $1$ is served and the  probabilities are $\delta/(\delta + (1-\pi') \lambda)$ and $(1-\pi') \lambda/(\delta + (1-\pi') \lambda)$ for a arrival and a catastrophe respectively, where 
$$\pi' := \mathbb{P}(\min_{0<t<\tilde{D}_1}Q(t)=0|\tilde{T}_1=0)$$ 
is the probability that the queue is empty before the first catastrophe occurs.
If this ``good'' event is a catastrophe, then $Q(\tilde{D}_1) =1$. If the first good event is an arrival, then customer 1 will not be served any more before $\tilde{D}_1$, and because $\tilde{D}_1$ is exponentially distributed, the number of individuals that arrived strictly after $\tilde{T}_1$ and will add to $Q(\tilde{D}_1)$, is distributed as $Q(\tilde{D}_1)$. So, 
$$\mathbb{P}(Q(\tilde{D}_1)=1|\tilde{T}_1=0 ,\min_{0<t<\tilde{D}_1}Q(t)>0) = \frac{\delta}{\delta +(1-\pi') \lambda}$$
and for $k \in \mathbb{N}$,
$$\mathbb{P}(Q(\tilde{D}_1)=k+1|\tilde{T}_1=0 ,\min_{0<t<\tilde{D}_1}Q(t)>0) = \frac{(1-\pi') \lambda}{\delta +(1-\pi')  \lambda}\mathbb{P}(Q(\tilde{D}_1)=k|\tilde{T}_1=0 ,\min_{0<t<\tilde{D}_1}Q(t)>0).$$
Combining this with the result of equation (\ref{relationtwoqueueing}) in section \ref{timegeneral} gives that for $k \in \mathbb{N}$, $\mathbb{P}(Q(\tilde{D}_1)=k|Q(0)=0,Q(\tilde{D}_1)>0) = p' (1-p')^{k-1}$, where
\begin{equation}\label{plifoform}
p' = \delta/(\delta + (1-\pi') \lambda).
\end{equation}
The same argument can be used to show that $\mathbb{P}(Q(\tilde{D}_1)=0|Q(0)=0) = p'$. 
Note that $\pi'$ is the probability that during the period customer 1 is served,  neither catastrophes nor arrivals of customers that are still in the queue at $\tilde{D}_1$, occur. Catastrophes occur at rate $\delta$ and arrivals of customers that are still in the queue after an exponentially ($\delta$) distributed time occur at rate $(1-\pi') \lambda$. So, $\pi'$ is the smallest (and unique) root in $[0,1]$ of 
\begin{equation}\label{pilifo}
\pi' = \mathbb{E}(e^{-(\delta + (1-\pi')\lambda)L}).
\end{equation}

This is the same equation as the equation for $\pi$ given in the previous section and thus, 
$$\mathbb{P}(Q^{\mbox{PS}}(\tilde{D}_1)=k|Q^{\mbox{PS}}(0)=0) = \mathbb{P}(Q^{\mbox{LIFO}}(\tilde{D}_1)=k|Q^{\mbox{LIFO}}(0)=0).$$ 
In particular, we have provided a straightforward proof for the fact that $\mathbb{P}(Q(\tilde{D}_1)=k|Q(0)=0,Q(\tilde{D}_1)>0) = p (1-p)^k$, in the $M/M/1$ queue with catastrophes and $\mathbb{P}(Z(D_1)=k|D_1 < \infty) = p (1-p)^k$ in the simple birth-death process. 

\section{Possible applications}\label{secapli}

\subsection{Within herd spread of contagious animal diseases}\label{sectionwithinheard}
The spread of very contagious animal diseases like Classical Swine Fever (CSF), Foot and Mouth Disease (FMD) and Avian Influenza (AI) within farms usually progresses as follows: The pathogen is imported in a farm and one animal in the farm becomes infected. Then an $SIR$ epidemic starts to spread within the farm. At a certain moment the farmer observes that an animal is diseased. Upon this detection the whole herd will be culled to prevent spread to other farms. In general the number of animals within a farm is large and the number of infectious animals at the moment of detection is small enough to accept the branching process approximation.

In \cite{Trap04} the following model for the spread of CSF within a farm is used. The spread starts with one infectious individual at time $t=0$. Every infected animal ``brings forth'' new infected  individuals at rate $\lambda$. Infectious animals recover at rate $\mu$ and are detected at rate $\delta$. Upon detection all animals in the farm are culled. Because the animals are immediately culled upon detection, the distribution of the number of infectious individuals at the moment of detection is the best information we can hope for to obtain if we only consider the within-farm spread. However, questions on this distribution are not addressed in \cite{Trap04}.

After some straightforward algebra (\ref{pandpi}) leads to 
$$\pi  = \frac{\mu + \delta + \lambda - \sqrt{(\mu + \delta + \lambda)^2 - 4 \lambda \mu}}{2 \lambda}$$
and
$$p = 1- \frac{\mu + \delta + \lambda - \sqrt{(\mu + \delta + \lambda)^2 - 4 \lambda \mu}}{2 \mu}.$$

Because the distribution of the number of infectious animals at the moment of detection is described by 1 parameter, more information than the number of infectious animals in a detected herd is needed in order to estimate $r := \lambda -\mu$ and $R := \lambda/\mu$ (respectively the exponential growth rate of the expected number of infected individuals and the basic reproduction number \cite{Diek00} within a farm). 
Possible further information should be obtained by contact tracing (by which we may obtain estimates for the real time between the first infection within a farm and the moment of detection), or by looking for traces of immune response in all culled animals on a farm and obtain the total number of animals infected before detection. Note that this total number is geometrically distributed with parameter $\delta/(\lambda + \delta)$, because every time a detection, infection or recovery occurs, the probability that a detection occurs is $\delta/\lambda$ times as high as the probability that an infection occurs.

\subsection{Spread of nosocomial pathogens like MRSA}

The spread of infectious diseases in hospitals is usually different from the spread of diseases in the general community.  Because of antibiotic pressure and because of the weak immune responses of hospitalized patients, these people are more susceptible to many pathogens, than healthy people. 
In \cite{Bootip} the following model for the spread of MRSA in hospitals in countries with low prevalence of MRSA (like The Nordic countries and The Netherlands) is used: At very low rate MRSA carriers enter a hospital. Say that the MRSA is brought in at time $t=0$ and carriers arrive as single patients, so $I(0) =1$. Because the rate at which MRSA carriers enter the hospital is very low, we exclude further introductions from outside the hospital during the outbreak started by the first entering of an MRSA carrier. Every infected patient infects other patients at rate $\lambda$, is discharged from the hospital at rate $\mu$ and is detected at rate $\delta$. As long as $I(t)$ is small compared to the total number of patients in the hospital or ward, the constant infection rate is reasonable.

Upon detection all patients in the hospital (or ward) are screened for MRSA colonisation and all infectious patients still in the hospital will be detected and isolated. Ideally this would make further infections within the hospital impossible. So $I(t)$ increases by 1 at rate $\lambda I(t)$, decreases by 1 at rate $\mu I(t)$ and detection occurs at rate $\delta I(t)$. The dynamics of $I(t)$ in this model are exactly the same as the dynamics of $I(t)$ in the model for the within farm spread of contagious animal diseases. Therefore, all results of the previous subsection \ref{sectionwithinheard} can be used for the spread of low prevalence nosocomial pathogens. Note that the assumption that patients leave the hospital at a constant rate is not necessary for applying the results of this paper.

\subsection{Change of behaviour because of knowledge of the epidemic}

If an infectious disease is known to spread in a certain region, physicians will be more alert on symptoms of the disease and people will  avoid crowded places or try to prevent their selves to become infected in other ways, like wearing masks as people did during the SARS epidemic in Asia. Therefore, it is reasonable to assume that the first detection of an infected person will lead to an increased rate of detection and a decreased rate of infection.   

Assume that the infection rate before the first detection was $\lambda_1 I(t)$ and that the detection rate was $\delta_1 I(t)$, while after the first detection the infection (resp. detection) rate will be $\lambda_2 I(t)$ (resp. $\delta_2 I(t)$).  We assume that the recovery rate per individual, $\mu$ does not change because of the first detection. From previous subsections we know that the number of infectious individuals at the moment of the first detection is geometrically distributed with parameter $p_1 := 1- \frac{\mu + \delta_1 + \lambda_1 - \sqrt{(\mu + \delta_1 + \lambda_1)^2 - 4 \lambda_1 \mu}}{2 \mu}$. 

From \cite[p.252]{Grim92} we deduce that the distribution of the number of individuals in a birth-death process at time $\tau$, $I(\tau)$, with per capita birth rate $\lambda_2$ and death rate $\mu$, which started with 1 infectious individual is given by $\mathbb{P}(Z(\tau)=k) = q_k(\tau)$ for $k \in \mathbb{N}_0$, with 

\begin{eqnarray*}
q_0(\tau) & = & \frac{e^{r_2 \tau}-1}{R_2 e^{r_2 \tau}-1}, \\
q_i(\tau) & = & (1-q_0(\tau)) (1-R_2 q_0(\tau))(R_2 q_0(\tau))^{i-1} \qquad
\mbox{for all} \ i \geq 1,
\end{eqnarray*}
where $r_2= \lambda_2 - \mu$ and $R_2= \frac{\lambda_2}{\mu}$.

Some algebra yields that if the initial number of infectious individuals is  geometrically distributed with parameter $p_1$, then at time $\tau$ the number of infectious individuals $Y(\tau)$ will be  given by:
\begin{eqnarray*}
\mathbb{P}(Y(\tau) =0) & =  & \frac{q_0(\tau) p_1}{1-(1-p_1)q_0(\tau)}\\
\mathbb{P}(Y(\tau) =i) & = &  (1-  \frac{q_0(\tau) p_1}{1-(1-p_1)q_0(\tau)}) \frac{(1-R_2 q_0(\tau)) p_1}{1-(1-p_1)q_0(\tau)} (1-  \frac{(1-R_2 q_0(\tau)) p_1}{1-(1-p_1)q_0(\tau)})^{i-1} \qquad
\mbox{for} \ i \geq 1
\end{eqnarray*}

In this example an exponentially distributed infectious period is assumed, because then the number of infectious individuals at the time of first detection contains as much information for the description of the progress of the epidemic after this detection, as knowledge of the whole process up to the time of first detection does. Note that for the epidemic with general infectious period it is possible to obtain the distribution of the number of infectious individuals that has been infectious for at least $x$ time units  at the moment of first detection from \cite[eq.(2.4)]{Kita93}.

\section{Extensions, limitations and questions}

In the previous sections we have assumed that the infection rate and detection rate during an infectious period are constant. It is tempting to conjecture that if the detection rate and infection rate are changing over time, but stay proportional, the number of infectious individuals at time $D_1$ is still geometrically distributed. Or formulated in the terminology of the corresponding branching process: let $\lambda(a)$ be the (possibly random) rate at which an individual at age $a$ gives birth and $\delta(a)$ the (possibly random) detection rate of this individual at age $a$ (if it is still alive at age $a$), where $\lambda(a)=c \delta(a)$ for some non-random constant $c$. 

An example of such a model is the $SEIR$ (Susceptible $\to$ Exposed $\to$ Infectious $\to$ Removed) epidemic model. In this model an individual first goes through a latent/exposed state after being infected, and after some random time, the individual becomes infectious itself. During the infectious period an individual can be detected, which happens at rate $\delta$. Apart from the latent period the model is the same as the $SIR$ epidemic model. 

If the infectious period is exponentially distributed with parameter $\mu$, and the latent periods are i.i.d.\ and distributed as $\Lambda$, then in the large population limit the dynamics of $E(t) + I(t)$ are described by: 
\begin{equation}\label{SEIReq}
\begin{array}{rcl}
\mathbb{P}(E(t+h)+I(t+h)= k+1|E(t)+I(t)=k,I(t)=l) & = & \lambda l h +o(h),\\ 
\mathbb{P}(E(t+h)+I(t+h)= k-1|E(t)+I(t)=k,I(t)=l) & = & \mu l h +o(h),\\
\mathbb{P}(D_i \in (t,t+h)|D_{i-1} \leq t < D_i, I(t)=l) & = & \delta l h + o(h),\\
\mathbb{P}(\mbox{more than 1 event in $(t,t+h)$}|I(t) = l) & = & o(h),
\end{array} 
\end{equation}
where events are infections, detections and recoveries. 
The time-change argument of Section \ref{secrela} can be applied with  
$$\tau'(t) =  \int_0^t \ind(I(t')>0)/I(t') dt'.$$
and we see that the dynamics of $I(\tau'(t)) + E(\tau'(t))$ are exactly the dynamics of an $M/M/1$ queue. 

However, simulations suggest that for the general $SEIR$ epidemic, with non-exponentially distributed infectious periods, $E(D_1)+I(D_1)$ is not geometrically distributed in the large population limit. 

An open question is whether there are $SEIR$ epidemics with non-exponentially distributed infectious periods, where $E(D_1)+I(D_1)$ is geometrically distributed in the large population limit? Or if $\lambda(a) = c \delta(a)$ is a deterministic function, for which functions $\lambda(a)$ and which distributions of the infectious period $L$, $I(D_1)$ is geometrically distributed? To answer these questions it might be helpful to have intuitive understanding of why 
\begin{equation}
\mathbb{P}(Q(\tilde{D}_1)=k|Q(0)=0,Q(\tilde{D}_1)>0) = p(1-p)^{k-1}
\end{equation}
holds for the general $M/G/1$-PS queue with catastrophes. The search for this intuitive understanding is still ongoing. 

\section*{Acknowledgements}
We thank O.~Boxma for helpful discussion and introducing us to the results obtained by Kitaev. We also thank A.~Sapozhnikov and O.~Diekmann for helpful discussion.
M.C.J.B. is supported by the Netherlands Organization for Scientific Research (VENI NWO Grant 916.86.128)


\end{document}